\documentstyle[12pt,amssymb,epic,eepic,fullpage]{amsart}\def\MRL{n}
\def\draft{n}
\theoremstyle{plain}

\newtheorem{theorem}{Theorem}

\newtheorem{proposition}{Proposition}[section]

\theoremstyle{definition}
\newtheorem{definition}[proposition]{Definition}

\theoremstyle{remark}

\newtheorem{remark}[proposition]{Remark}

\def\printname#1{
	\if\draft y
		\smash{\makebox[0pt]{\hspace{-0.5in}
			\raisebox{8pt}{\tt\tiny #1}}}
	\fi
}

\newcommand{\mathmode}[1]{$#1$}

\newlength{\standardunitlength}
\catcode`\@=11
\long\def\@makecaption#1#2{%
    \vskip 10pt
    \setbox\@tempboxa\hbox{
      \small\sf{\bfcaptionfont #1. }\ignorespaces #2}%
    \ifdim \wd\@tempboxa >\captionwidth {%
        \rightskip=\@captionmargin\leftskip=\@captionmargin
        \unhbox\@tempboxa\par}%
      \else
        \hbox to\hsize{\hfil\box\@tempboxa\hfil}%
    \fi}
\font\bfcaptionfont=cmssbx10 scaled \magstephalf
\newdimen\@captionmargin\@captionmargin=2\parindent
\newdimen\captionwidth\captionwidth=\hsize
\catcode`\@=12

\def\lbl#1{\label{#1}\printname{#1}}
\def\biblbl#1{\bibitem{#1}\printname{#1}}

\theoremstyle{plain}
\newtheorem{Lemma}{Lemma}
\newtheorem{SubLemma}{SubLemma}[Lemma]
\newtheorem{SubSubLemma}{SubSubLemma}[SubLemma]

\def\Acknowledgement{
  I would like to thank P.~Etingof, Y.~Karshon, X-S.~Lin, A.~Referee,
  and J.~Wunsch for their comments and suggestions, X-S.~Lin and
  Z.~Wang for suggesting the problem in their
  paper~\cite{LinWang:IntegralGeometry}, and D.~Zeilberger for his
  contributions to the theory of SubSubLemmas
  (\cite{Zeilberger:SubSubLemmas}). I am the littlest of his
  disciples.
}
\newcommand{\EEPIC}[2]{
	\setlength{\unitlength}{#2\standardunitlength}
        \begin{array}{c}  \hspace{-1.7mm}
                \raisebox{-8pt}{#1}
                \hspace{-1.9mm}
        \end{array}
}
\newcommand{\silentEEPIC}[2]
        {\setlength{\unitlength}{#2\standardunitlength}
        \begin{array}{c}  \hspace{-1.7mm}
                \raisebox{-2pt}{#1}
                \hspace{-1.9mm}
        \end{array}
}
\def\Cutting
{{
\begin{picture}(550,155)(0,-10)
\dashline{7.000}(145,15)(100,15)(75,15)(70,30)
\path(74.427,23.043)(70.000,30.000)(70.632,21.778)
\dottedline{5}(55,0)(55,20)(55,40)
	(55,80)(80,120)(120,120)
	(120,80)(90,60)(85,40)
	(85,20)(85,0)
\path(83.000,8.000)(85.000,0.000)(87.000,8.000)
\dashline{7.000}(215,60)(170,60)(125,75)
\path(133.222,74.368)(125.000,75.000)(131.957,70.573)
\dottedline{5}(55,30)(60,35)(70,35)
	(80,35)(85,30)
\path(0,100)(140,100)
\path(0,40)(140,40)
\path(40,140)(40,0)
\path(100,140)(100,0)
\dashline{7.000}(215,125)(160,125)(140,105)
\path(144.243,112.071)(140.000,105.000)(147.071,109.243)
\dashline{7.000}(325,125)(380,125)(410,105)
\path(402.234,107.774)(410.000,105.000)(404.453,111.102)
\dashline{7.000}(335,60)(460,60)
\path(452.000,58.000)(460.000,60.000)(452.000,62.000)
\dashline{7.000}(405,15)(460,15)(475,35)
\path(471.800,27.400)(475.000,35.000)(468.600,29.800)
\path(410,100)(550,100)
\path(410,40)(550,40)
\path(450,140)(450,0)
\path(510,140)(510,0)
\dottedline{5}(465,0)(465,80)(490,120)
	(530,120)(530,20)(525,15)
	(500,15)(495,20)(495,70)
\path(497.000,62.000)(495.000,70.000)(493.000,62.000)
\dottedline{5}(465,35)(495,45)
\put(225,55){\makebox(0,0)[lb]{
\raisebox{0pt}[0pt][0pt]{\shortstack[l]{{\rm\scriptsize The original
$\gamma_i$}}}}}
\put(155,10){\makebox(0,0)[lb]{
\raisebox{0pt}[0pt][0pt]{\shortstack[l]{{\rm\scriptsize A short cut with fewer
intersections}}}}}
\put(225,120){\makebox(0,0)[lb]{
\raisebox{0pt}[0pt][0pt]{\shortstack[l]{{\rm\scriptsize A part of $\pi K$}}}}}
\end{picture}
}}
\def\Moves
{{
\begin{picture}(240,95)(0,-10)
\dottedline{5}(0,80)(80,0)
\path(72.929,4.243)(80.000,0.000)(75.757,7.071)
\dottedline{5}(80,80)(0,0)
\path(4.243,7.071)(0.000,0.000)(7.071,4.243)
\dottedline{5}(160,80)(190,60)(190,20)(160,0)
\path(165.547,6.102)(160.000,0.000)(167.766,2.774)
\dottedline{5}(240,80)(210,60)(210,20)(240,0)
\path(232.234,2.774)(240.000,0.000)(234.453,6.102)
\path(100,40)(140,40)
\path(132.000,38.000)(140.000,40.000)(132.000,42.000)
\end{picture}
}}
\def\OddChoice
{{
\begin{picture}(500,161)(0,-10)
\put(360.000,20.500){\arc{40.311}{0.1244}{3.0172}}
\put(380.000,20.500){\arc{40.311}{0.1244}{3.0172}}
\put(140.000,20.500){\arc{40.311}{0.1244}{3.0172}}
\put(160.000,20.500){\arc{40.311}{0.1244}{3.0172}}
\put(370,3){\blacken\ellipse{6}{6}}
\put(370,3){\ellipse{6}{6}}
\put(150,3){\blacken\ellipse{6}{6}}
\put(150,3){\ellipse{6}{6}}
\put(270,108){\blacken\ellipse{6}{6}}
\put(270,108){\ellipse{6}{6}}
\put(270,3){\blacken\ellipse{6}{6}}
\put(270,3){\ellipse{6}{6}}
\dashline{7.000}(500,18)(0,18)
\path(280,118)(260,98)
\path(280,118)(260,98)
\path(260,118)(280,98)
\path(260,118)(280,98)
\dottedline{5}(270,108)(220,98)(210,78)
	(240,38)(270,3)
\path(491.924,21.334)(500.000,23.000)(492.090,25.331)
\path(500,23)	(495.820,23.184)
	(491.933,23.376)
	(488.323,23.579)
	(484.977,23.794)
	(481.878,24.024)
	(479.014,24.269)
	(476.368,24.532)
	(473.927,24.815)
	(469.600,25.446)
	(465.914,26.177)
	(462.753,27.024)
	(460.000,28.000)

\path(460,28)	(455.502,30.447)
	(450.642,34.064)
	(448.111,36.186)
	(445.530,38.452)
	(442.913,40.810)
	(440.274,43.211)
	(437.626,45.606)
	(434.984,47.945)
	(432.360,50.177)
	(429.768,52.253)
	(427.223,54.123)
	(424.737,55.738)
	(420.000,58.000)

\path(420,58)	(415.567,59.141)
	(412.927,59.571)
	(410.156,59.829)
	(407.364,59.859)
	(404.666,59.603)
	(400.000,58.000)

\path(400,58)	(396.854,55.719)
	(394.010,52.895)
	(391.409,49.411)
	(388.993,45.151)
	(387.835,42.693)
	(386.702,39.997)
	(385.586,37.048)
	(384.479,33.832)
	(383.374,30.334)
	(382.264,26.539)
	(381.142,22.432)
	(380.000,18.000)

\path(340,18)	(339.241,22.073)
	(338.549,25.579)
	(337.910,28.577)
	(337.310,31.126)
	(336.166,35.110)
	(335.000,38.000)

\path(335,38)	(333.725,40.428)
	(332.167,43.079)
	(330.360,45.918)
	(328.339,48.907)
	(326.138,52.008)
	(323.792,55.186)
	(321.337,58.402)
	(318.805,61.620)
	(316.233,64.804)
	(313.655,67.915)
	(311.106,70.918)
	(308.619,73.774)
	(306.231,76.448)
	(303.975,78.901)
	(300.000,83.000)

\path(300,83)	(297.234,85.482)
	(293.291,88.575)
	(290.732,90.460)
	(287.703,92.630)
	(284.145,95.129)
	(280.000,98.000)

\path(260,118)	(257.426,121.336)
	(255.140,124.166)
	(251.254,128.487)
	(247.991,131.315)
	(245.000,133.000)

\path(245,133)	(241.355,134.204)
	(237.096,135.096)
	(232.426,135.660)
	(227.546,135.880)
	(222.662,135.740)
	(217.974,135.224)
	(213.685,134.316)
	(210.000,133.000)

\path(210,133)	(205.577,130.365)
	(201.120,126.668)
	(196.742,122.191)
	(192.556,117.216)
	(190.571,114.631)
	(188.676,112.027)
	(186.886,109.440)
	(185.214,106.906)
	(182.284,102.136)
	(180.000,98.000)

\path(180,98)	(178.730,95.418)
	(177.497,92.651)
	(176.293,89.668)
	(175.112,86.442)
	(173.946,82.941)
	(172.787,79.138)
	(171.628,75.002)
	(170.462,70.505)
	(169.282,65.618)
	(168.684,63.018)
	(168.080,60.310)
	(167.468,57.489)
	(166.848,54.553)
	(166.220,51.496)
	(165.581,48.317)
	(164.931,45.010)
	(164.270,41.573)
	(163.595,38.001)
	(162.907,34.291)
	(162.205,30.440)
	(161.486,26.444)
	(160.752,22.298)
	(160.000,18.000)

\path(120,18)	(119.565,21.675)
	(119.419,24.411)
	(120.000,28.000)

\path(120,28)	(121.405,31.049)
	(123.485,34.321)
	(126.322,36.933)
	(130.000,38.000)

\path(130,38)	(133.678,36.933)
	(136.515,34.322)
	(138.595,31.049)
	(140.000,28.000)

\path(140,28)	(140.581,24.411)
	(140.435,21.675)
	(140.000,18.000)

\path(180,18)	(180.579,21.083)
	(181.142,23.733)
	(182.278,27.910)
	(183.525,30.882)
	(185.000,33.000)

\path(185,33)	(189.957,37.476)
	(192.968,39.606)
	(196.265,41.658)
	(199.801,43.624)
	(203.523,45.499)
	(207.382,47.278)
	(211.328,48.955)
	(215.311,50.524)
	(219.280,51.979)
	(223.185,53.315)
	(226.977,54.525)
	(230.604,55.605)
	(234.017,56.547)
	(237.166,57.348)
	(240.000,58.000)

\path(240,58)	(242.546,58.503)
	(245.250,58.939)
	(248.138,59.308)
	(251.239,59.609)
	(254.579,59.844)
	(258.187,60.012)
	(262.089,60.112)
	(266.313,60.146)
	(270.887,60.112)
	(275.838,60.012)
	(278.463,59.936)
	(281.193,59.844)
	(284.031,59.735)
	(286.981,59.609)
	(290.045,59.467)
	(293.228,59.308)
	(296.532,59.132)
	(299.962,58.939)
	(303.520,58.729)
	(307.210,58.503)
	(311.035,58.260)
	(315.000,58.000)

\path(325,58)	(328.320,56.332)
	(331.396,54.755)
	(334.238,53.260)
	(336.859,51.840)
	(341.479,49.196)
	(345.344,46.766)
	(348.540,44.489)
	(351.156,42.309)
	(355.000,38.000)

\path(355,38)	(356.581,35.102)
	(357.863,31.115)
	(358.429,28.567)
	(358.964,25.570)
	(359.483,22.068)
	(360.000,18.000)

\path(400,18)	(400.000,21.034)
	(400.000,23.648)
	(400.000,27.797)
	(400.000,30.797)
	(400.000,33.000)

\path(400,33)	(400.000,35.203)
	(400.000,38.203)
	(400.000,42.352)
	(400.000,44.966)
	(400.000,48.000)

\path(400,63)	(399.296,67.073)
	(398.645,70.579)
	(398.030,73.577)
	(397.438,76.126)
	(396.262,80.110)
	(395.000,83.000)

\path(395,83)	(393.052,86.466)
	(390.767,90.245)
	(388.173,94.280)
	(385.301,98.511)
	(382.181,102.878)
	(378.842,107.322)
	(375.314,111.784)
	(371.628,116.206)
	(367.813,120.527)
	(363.898,124.689)
	(359.914,128.633)
	(355.891,132.298)
	(351.858,135.627)
	(347.845,138.560)
	(343.883,141.037)
	(340.000,143.000)

\path(340,143)	(336.330,144.258)
	(332.049,145.139)
	(327.362,145.652)
	(322.474,145.805)
	(317.589,145.607)
	(312.911,145.068)
	(308.647,144.196)
	(305.000,143.000)

\path(305,143)	(300.014,140.191)
	(297.388,138.110)
	(294.576,135.479)
	(291.504,132.226)
	(288.100,128.277)
	(284.290,123.560)
	(282.209,120.890)
	(280.000,118.000)

\path(260,98)	(256.526,94.022)
	(253.857,91.185)
	(250.000,88.000)

\path(250,88)	(245.180,85.595)
	(242.308,84.433)
	(239.191,83.301)
	(235.878,82.202)
	(232.415,81.138)
	(228.849,80.113)
	(225.229,79.128)
	(221.602,78.185)
	(218.015,77.288)
	(214.516,76.439)
	(211.152,75.640)
	(207.971,74.895)
	(205.021,74.204)
	(202.348,73.572)
	(200.000,73.000)

\path(200,73)	(196.347,72.152)
	(191.353,71.113)
	(188.170,70.487)
	(184.433,69.769)
	(180.067,68.945)
	(177.626,68.488)
	(175.000,68.000)

\path(160,63)	(155.812,62.754)
	(151.918,62.508)
	(148.303,62.260)
	(144.951,62.008)
	(141.849,61.749)
	(138.982,61.484)
	(136.335,61.208)
	(133.894,60.922)
	(129.568,60.307)
	(125.889,59.625)
	(122.739,58.861)
	(120.000,58.000)

\path(120,58)	(115.648,55.990)
	(110.741,53.024)
	(108.137,51.283)
	(105.463,49.422)
	(102.743,47.481)
	(100.000,45.500)
	(97.257,43.519)
	(94.537,41.578)
	(91.863,39.717)
	(89.259,37.976)
	(84.352,35.010)
	(80.000,33.000)

\path(80,33)	(77.350,32.113)
	(74.523,31.279)
	(71.491,30.492)
	(68.223,29.751)
	(64.691,29.050)
	(60.864,28.386)
	(56.715,27.757)
	(52.214,27.157)
	(47.331,26.584)
	(44.737,26.306)
	(42.036,26.033)
	(39.226,25.765)
	(36.302,25.501)
	(33.260,25.241)
	(30.098,24.985)
	(26.811,24.731)
	(23.395,24.480)
	(19.847,24.231)
	(16.164,23.984)
	(12.341,23.737)
	(8.375,23.492)
	(4.263,23.246)
	(0.000,23.000)

\end{picture}
}}
\def\Pulling
{{
\begin{picture}(500,168)(0,-10)
\put(360.000,27.500){\arc{40.311}{0.1244}{3.0172}}
\put(380.000,27.500){\arc{40.311}{0.1244}{3.0172}}
\put(140.000,27.500){\arc{40.311}{0.1244}{3.0172}}
\put(160.000,27.500){\arc{40.311}{0.1244}{3.0172}}
\put(370,10){\blacken\ellipse{6}{6}}
\put(370,10){\ellipse{6}{6}}
\put(150,10){\blacken\ellipse{6}{6}}
\put(150,10){\ellipse{6}{6}}
\put(270,10){\blacken\ellipse{6}{6}}
\put(270,10){\ellipse{6}{6}}
\dashline{7.000}(500,25)(0,25)
\path(260,105)(230,100)(220,85)
	(275,15)(270,10)
\path(260,125)(220,110)(205,85)
	(262,1)(270,0)(275,5)(270,10)
\path(280,125)(275,120)(223,107)
	(210,85)(265,5)(270,10)
\path(280,105)(270,113)(227,103)
	(215,85)(265,15)(270,10)
\path(491.924,28.334)(500.000,30.000)(492.090,32.331)
\path(500,30)	(495.820,30.184)
	(491.933,30.376)
	(488.323,30.579)
	(484.977,30.794)
	(481.878,31.024)
	(479.014,31.269)
	(476.368,31.532)
	(473.927,31.815)
	(469.600,32.446)
	(465.914,33.177)
	(462.753,34.024)
	(460.000,35.000)

\path(460,35)	(455.502,37.447)
	(450.642,41.064)
	(448.111,43.186)
	(445.530,45.452)
	(442.913,47.810)
	(440.274,50.211)
	(437.626,52.606)
	(434.984,54.945)
	(432.360,57.177)
	(429.768,59.253)
	(427.223,61.123)
	(424.737,62.738)
	(420.000,65.000)

\path(420,65)	(415.567,66.141)
	(412.927,66.571)
	(410.156,66.829)
	(407.364,66.859)
	(404.666,66.603)
	(400.000,65.000)

\path(400,65)	(396.854,62.719)
	(394.010,59.895)
	(391.409,56.411)
	(388.993,52.151)
	(387.835,49.693)
	(386.702,46.997)
	(385.586,44.048)
	(384.479,40.832)
	(383.374,37.334)
	(382.264,33.539)
	(381.142,29.432)
	(380.000,25.000)

\path(340,25)	(339.241,29.073)
	(338.549,32.579)
	(337.910,35.577)
	(337.310,38.126)
	(336.166,42.110)
	(335.000,45.000)

\path(335,45)	(333.725,47.428)
	(332.167,50.079)
	(330.360,52.918)
	(328.339,55.907)
	(326.138,59.008)
	(323.792,62.186)
	(321.337,65.402)
	(318.805,68.620)
	(316.233,71.804)
	(313.655,74.915)
	(311.106,77.918)
	(308.619,80.774)
	(306.231,83.448)
	(303.975,85.901)
	(300.000,90.000)

\path(300,90)	(297.234,92.482)
	(293.291,95.575)
	(290.732,97.460)
	(287.703,99.630)
	(284.145,102.129)
	(280.000,105.000)

\path(260,125)	(257.426,128.336)
	(255.140,131.166)
	(251.254,135.487)
	(247.991,138.315)
	(245.000,140.000)

\path(245,140)	(241.355,141.204)
	(237.096,142.096)
	(232.426,142.660)
	(227.546,142.880)
	(222.662,142.740)
	(217.974,142.224)
	(213.685,141.316)
	(210.000,140.000)

\path(210,140)	(205.577,137.365)
	(201.120,133.668)
	(196.742,129.191)
	(192.556,124.216)
	(190.571,121.631)
	(188.676,119.027)
	(186.886,116.440)
	(185.214,113.906)
	(182.284,109.136)
	(180.000,105.000)

\path(180,105)	(178.730,102.418)
	(177.497,99.651)
	(176.293,96.668)
	(175.112,93.442)
	(173.946,89.941)
	(172.787,86.138)
	(171.628,82.002)
	(170.462,77.505)
	(169.282,72.618)
	(168.684,70.018)
	(168.080,67.310)
	(167.468,64.489)
	(166.848,61.553)
	(166.220,58.496)
	(165.581,55.317)
	(164.931,52.010)
	(164.270,48.573)
	(163.595,45.001)
	(162.907,41.291)
	(162.205,37.440)
	(161.486,33.444)
	(160.752,29.298)
	(160.000,25.000)

\path(120,25)	(119.565,28.675)
	(119.419,31.411)
	(120.000,35.000)

\path(120,35)	(121.405,38.049)
	(123.485,41.321)
	(126.322,43.933)
	(130.000,45.000)

\path(130,45)	(133.678,43.933)
	(136.515,41.322)
	(138.595,38.049)
	(140.000,35.000)

\path(140,35)	(140.581,31.411)
	(140.435,28.675)
	(140.000,25.000)

\path(180,25)	(180.579,28.083)
	(181.142,30.733)
	(182.278,34.910)
	(183.525,37.882)
	(185.000,40.000)

\path(185,40)	(189.957,44.476)
	(192.968,46.606)
	(196.265,48.658)
	(199.801,50.624)
	(203.523,52.499)
	(207.382,54.278)
	(211.328,55.955)
	(215.311,57.524)
	(219.280,58.979)
	(223.185,60.315)
	(226.977,61.525)
	(230.604,62.605)
	(234.017,63.547)
	(237.166,64.348)
	(240.000,65.000)

\path(240,65)	(242.546,65.503)
	(245.250,65.939)
	(248.138,66.308)
	(251.239,66.609)
	(254.579,66.844)
	(258.187,67.012)
	(262.089,67.112)
	(266.313,67.146)
	(270.887,67.112)
	(275.838,67.012)
	(278.463,66.936)
	(281.193,66.844)
	(284.031,66.735)
	(286.981,66.609)
	(290.045,66.467)
	(293.228,66.308)
	(296.532,66.132)
	(299.962,65.939)
	(303.520,65.729)
	(307.210,65.503)
	(311.035,65.260)
	(315.000,65.000)

\path(325,65)	(328.320,63.332)
	(331.396,61.755)
	(334.238,60.260)
	(336.859,58.840)
	(341.479,56.196)
	(345.344,53.766)
	(348.540,51.489)
	(351.156,49.309)
	(355.000,45.000)

\path(355,45)	(356.581,42.102)
	(357.863,38.115)
	(358.429,35.567)
	(358.964,32.570)
	(359.483,29.068)
	(360.000,25.000)

\path(400,25)	(400.000,28.034)
	(400.000,30.648)
	(400.000,34.797)
	(400.000,37.797)
	(400.000,40.000)

\path(400,40)	(400.000,42.203)
	(400.000,45.203)
	(400.000,49.352)
	(400.000,51.966)
	(400.000,55.000)

\path(400,70)	(399.296,74.073)
	(398.645,77.579)
	(398.030,80.577)
	(397.438,83.126)
	(396.262,87.110)
	(395.000,90.000)

\path(395,90)	(393.052,93.466)
	(390.767,97.245)
	(388.173,101.280)
	(385.301,105.511)
	(382.181,109.878)
	(378.842,114.322)
	(375.314,118.784)
	(371.628,123.206)
	(367.813,127.527)
	(363.898,131.689)
	(359.914,135.633)
	(355.891,139.298)
	(351.858,142.627)
	(347.845,145.560)
	(343.883,148.037)
	(340.000,150.000)

\path(340,150)	(336.330,151.258)
	(332.049,152.139)
	(327.362,152.652)
	(322.474,152.805)
	(317.589,152.607)
	(312.911,152.068)
	(308.647,151.196)
	(305.000,150.000)

\path(305,150)	(300.014,147.191)
	(297.388,145.110)
	(294.576,142.479)
	(291.504,139.226)
	(288.100,135.277)
	(284.290,130.560)
	(282.209,127.890)
	(280.000,125.000)

\path(260,105)	(256.526,101.022)
	(253.857,98.185)
	(250.000,95.000)

\path(250,95)	(245.180,92.595)
	(242.308,91.433)
	(239.191,90.301)
	(235.878,89.202)
	(232.415,88.138)
	(228.849,87.113)
	(225.229,86.128)
	(221.602,85.185)
	(218.015,84.288)
	(214.516,83.439)
	(211.152,82.640)
	(207.971,81.895)
	(205.021,81.204)
	(202.348,80.572)
	(200.000,80.000)

\path(200,80)	(196.347,79.152)
	(191.353,78.113)
	(188.170,77.487)
	(184.433,76.769)
	(180.067,75.945)
	(177.626,75.488)
	(175.000,75.000)

\path(160,70)	(155.812,69.754)
	(151.918,69.508)
	(148.303,69.260)
	(144.951,69.008)
	(141.849,68.749)
	(138.982,68.484)
	(136.335,68.208)
	(133.894,67.922)
	(129.568,67.307)
	(125.889,66.625)
	(122.739,65.861)
	(120.000,65.000)

\path(120,65)	(115.648,62.990)
	(110.741,60.024)
	(108.137,58.283)
	(105.463,56.422)
	(102.743,54.481)
	(100.000,52.500)
	(97.257,50.519)
	(94.537,48.578)
	(91.863,46.717)
	(89.259,44.976)
	(84.352,42.010)
	(80.000,40.000)

\path(80,40)	(77.350,39.113)
	(74.523,38.279)
	(71.491,37.492)
	(68.223,36.751)
	(64.691,36.050)
	(60.864,35.386)
	(56.715,34.757)
	(52.214,34.157)
	(47.331,33.584)
	(44.737,33.306)
	(42.036,33.033)
	(39.226,32.765)
	(36.302,32.501)
	(33.260,32.241)
	(30.098,31.985)
	(26.811,31.731)
	(23.395,31.480)
	(19.847,31.231)
	(16.164,30.984)
	(12.341,30.737)
	(8.375,30.492)
	(4.263,30.246)
	(0.000,30.000)

\end{picture}
}}
\def\SideView
{{
\begin{picture}(280,295)(0,-10)
\path(165,165)(195,195)
\path(205,205)(235,235)
\path(245,245)(280,280)
\path(0,0)(40,40)(80,80)
	(120,120)(155,155)
\dashline{7.000}(0,0)(280,0)
\path(272.000,-2.000)(280.000,0.000)(272.000,2.000)
\dashline{7.000}(0,0)(0,280)
\path(2.000,272.000)(0.000,280.000)(-2.000,272.000)
\path(155,155)	(155.351,150.276)
	(155.689,145.950)
	(156.014,142.009)
	(156.330,138.441)
	(156.637,135.232)
	(156.938,132.370)
	(157.529,127.635)
	(158.115,124.132)
	(158.713,121.761)
	(160.000,120.000)

\path(160,120)	(161.577,122.306)
	(162.247,125.365)
	(162.858,129.872)
	(163.145,132.708)
	(163.424,135.956)
	(163.695,139.632)
	(163.960,143.751)
	(164.221,148.331)
	(164.351,150.799)
	(164.480,153.389)
	(164.610,156.101)
	(164.739,158.939)
	(164.869,161.905)
	(165.000,165.000)

\path(195,195)	(195.041,192.290)
	(195.082,189.607)
	(195.122,186.953)
	(195.163,184.327)
	(195.203,181.729)
	(195.243,179.159)
	(195.283,176.616)
	(195.323,174.102)
	(195.402,169.154)
	(195.480,164.316)
	(195.558,159.586)
	(195.635,154.963)
	(195.712,150.446)
	(195.788,146.035)
	(195.864,141.729)
	(195.940,137.527)
	(196.015,133.427)
	(196.090,129.430)
	(196.164,125.533)
	(196.238,121.737)
	(196.312,118.040)
	(196.386,114.442)
	(196.459,110.941)
	(196.532,107.537)
	(196.605,104.228)
	(196.678,101.015)
	(196.751,97.896)
	(196.824,94.869)
	(196.896,91.935)
	(196.969,89.092)
	(197.042,86.340)
	(197.114,83.677)
	(197.187,81.103)
	(197.260,78.617)
	(197.406,73.905)
	(197.553,69.533)
	(197.701,65.494)
	(197.850,61.782)
	(198.000,58.390)
	(198.153,55.309)
	(198.307,52.534)
	(198.621,47.870)
	(198.946,44.341)
	(199.283,41.891)
	(200.000,40.000)

\path(200,40)	(200.760,42.015)
	(201.113,44.623)
	(201.451,48.380)
	(201.777,53.345)
	(201.935,56.300)
	(202.090,59.579)
	(202.243,63.191)
	(202.394,67.142)
	(202.543,71.441)
	(202.691,76.095)
	(202.764,78.558)
	(202.836,81.112)
	(202.909,83.758)
	(202.981,86.498)
	(203.053,89.333)
	(203.124,92.263)
	(203.196,95.289)
	(203.267,98.412)
	(203.338,101.634)
	(203.409,104.954)
	(203.480,108.375)
	(203.551,111.897)
	(203.622,115.520)
	(203.693,119.247)
	(203.763,123.078)
	(203.834,127.013)
	(203.905,131.054)
	(203.976,135.202)
	(204.048,139.457)
	(204.119,143.821)
	(204.191,148.294)
	(204.263,152.878)
	(204.335,157.573)
	(204.407,162.381)
	(204.480,167.302)
	(204.516,169.805)
	(204.553,172.337)
	(204.590,174.898)
	(204.626,177.487)
	(204.663,180.106)
	(204.700,182.754)
	(204.737,185.431)
	(204.774,188.137)
	(204.812,190.873)
	(204.849,193.639)
	(204.887,196.434)
	(204.924,199.260)
	(204.962,202.115)
	(205.000,205.000)

\path(235,235)	(235.041,232.290)
	(235.082,229.607)
	(235.122,226.953)
	(235.163,224.327)
	(235.203,221.729)
	(235.243,219.159)
	(235.283,216.616)
	(235.323,214.102)
	(235.402,209.154)
	(235.480,204.316)
	(235.558,199.586)
	(235.635,194.963)
	(235.712,190.446)
	(235.788,186.035)
	(235.864,181.729)
	(235.940,177.527)
	(236.015,173.427)
	(236.090,169.430)
	(236.164,165.533)
	(236.238,161.737)
	(236.312,158.040)
	(236.386,154.442)
	(236.459,150.941)
	(236.532,147.537)
	(236.605,144.228)
	(236.678,141.015)
	(236.751,137.896)
	(236.824,134.869)
	(236.896,131.935)
	(236.969,129.092)
	(237.042,126.340)
	(237.114,123.677)
	(237.187,121.103)
	(237.260,118.617)
	(237.406,113.905)
	(237.553,109.533)
	(237.701,105.494)
	(237.850,101.782)
	(238.000,98.390)
	(238.153,95.309)
	(238.307,92.534)
	(238.621,87.870)
	(238.946,84.341)
	(239.283,81.891)
	(240.000,80.000)

\path(240,80)	(240.760,82.015)
	(241.113,84.623)
	(241.451,88.380)
	(241.777,93.345)
	(241.935,96.300)
	(242.090,99.579)
	(242.243,103.191)
	(242.394,107.142)
	(242.543,111.441)
	(242.691,116.095)
	(242.764,118.558)
	(242.836,121.112)
	(242.909,123.758)
	(242.981,126.498)
	(243.053,129.333)
	(243.124,132.263)
	(243.196,135.289)
	(243.267,138.412)
	(243.338,141.634)
	(243.409,144.954)
	(243.480,148.375)
	(243.551,151.897)
	(243.622,155.520)
	(243.693,159.247)
	(243.763,163.078)
	(243.834,167.013)
	(243.905,171.054)
	(243.976,175.202)
	(244.048,179.457)
	(244.119,183.821)
	(244.191,188.294)
	(244.263,192.878)
	(244.335,197.573)
	(244.407,202.381)
	(244.480,207.302)
	(244.516,209.805)
	(244.553,212.337)
	(244.590,214.898)
	(244.626,217.487)
	(244.663,220.106)
	(244.700,222.754)
	(244.737,225.431)
	(244.774,228.137)
	(244.812,230.873)
	(244.849,233.639)
	(244.887,236.434)
	(244.924,239.260)
	(244.962,242.115)
	(245.000,245.000)

\put(260,5){\makebox(0,0)[lb]{
\raisebox{0pt}[0pt][0pt]{\shortstack[l]{\mathmode{s}}}}}
\put(5,260){\makebox(0,0)[lb]{
\raisebox{0pt}[0pt][0pt]{\shortstack[l]{\mathmode{z}}}}}
\end{picture}
}}
\def\VertView
{{
\begin{picture}(640,163)(0,-10)
\put(360.000,22.500){\arc{40.311}{0.1244}{3.0172}}
\put(380.000,22.500){\arc{40.311}{0.1244}{3.0172}}
\put(140.000,22.500){\arc{40.311}{0.1244}{3.0172}}
\put(160.000,22.500){\arc{40.311}{0.1244}{3.0172}}
\put(570.000,80.000){\arc{20.000}{3.1416}{6.2832}}
\put(560.000,80.000){\arc{80.000}{3.1416}{6.2832}}
\put(580.000,80.000){\arc{80.000}{3.1416}{6.2832}}
\put(370,5){\blacken\ellipse{6}{6}}
\put(370,5){\ellipse{6}{6}}
\put(150,5){\blacken\ellipse{6}{6}}
\put(150,5){\ellipse{6}{6}}
\put(270,110){\blacken\ellipse{6}{6}}
\put(270,110){\ellipse{6}{6}}
\dashline{7.000}(500,20)(0,20)
\path(280,120)(260,100)
\path(280,120)(260,100)
\path(260,120)(280,100)
\path(260,120)(280,100)
\path(500,80)(640,80)
\path(632.000,78.000)(640.000,80.000)(632.000,82.000)
\path(491.924,23.334)(500.000,25.000)(492.090,27.331)
\path(500,25)	(495.820,25.184)
	(491.933,25.376)
	(488.323,25.579)
	(484.977,25.794)
	(481.878,26.024)
	(479.014,26.269)
	(476.368,26.532)
	(473.927,26.815)
	(469.600,27.446)
	(465.914,28.177)
	(462.753,29.024)
	(460.000,30.000)

\path(460,30)	(455.502,32.447)
	(450.642,36.064)
	(448.111,38.186)
	(445.530,40.452)
	(442.913,42.810)
	(440.274,45.211)
	(437.626,47.606)
	(434.984,49.945)
	(432.360,52.177)
	(429.768,54.253)
	(427.223,56.123)
	(424.737,57.738)
	(420.000,60.000)

\path(420,60)	(415.567,61.141)
	(412.927,61.571)
	(410.156,61.829)
	(407.364,61.859)
	(404.666,61.603)
	(400.000,60.000)

\path(400,60)	(396.854,57.719)
	(394.010,54.895)
	(391.409,51.411)
	(388.993,47.151)
	(387.835,44.693)
	(386.702,41.997)
	(385.586,39.048)
	(384.479,35.832)
	(383.374,32.334)
	(382.264,28.539)
	(381.142,24.432)
	(380.000,20.000)

\path(340,20)	(339.241,24.073)
	(338.549,27.579)
	(337.910,30.577)
	(337.310,33.126)
	(336.166,37.110)
	(335.000,40.000)

\path(335,40)	(333.725,42.428)
	(332.167,45.079)
	(330.360,47.918)
	(328.339,50.907)
	(326.138,54.008)
	(323.792,57.186)
	(321.337,60.402)
	(318.805,63.620)
	(316.233,66.804)
	(313.655,69.915)
	(311.106,72.918)
	(308.619,75.774)
	(306.231,78.448)
	(303.975,80.901)
	(300.000,85.000)

\path(300,85)	(297.234,87.482)
	(293.291,90.575)
	(290.732,92.460)
	(287.703,94.630)
	(284.145,97.129)
	(280.000,100.000)

\path(260,120)	(257.426,123.336)
	(255.140,126.166)
	(251.254,130.487)
	(247.991,133.315)
	(245.000,135.000)

\path(245,135)	(241.355,136.204)
	(237.096,137.096)
	(232.426,137.660)
	(227.546,137.880)
	(222.662,137.740)
	(217.974,137.224)
	(213.685,136.316)
	(210.000,135.000)

\path(210,135)	(205.577,132.365)
	(201.120,128.668)
	(196.742,124.191)
	(192.556,119.216)
	(190.571,116.631)
	(188.676,114.027)
	(186.886,111.440)
	(185.214,108.906)
	(182.284,104.136)
	(180.000,100.000)

\path(180,100)	(178.730,97.418)
	(177.497,94.651)
	(176.293,91.668)
	(175.112,88.442)
	(173.946,84.941)
	(172.787,81.138)
	(171.628,77.002)
	(170.462,72.505)
	(169.282,67.618)
	(168.684,65.018)
	(168.080,62.310)
	(167.468,59.489)
	(166.848,56.553)
	(166.220,53.496)
	(165.581,50.317)
	(164.931,47.010)
	(164.270,43.573)
	(163.595,40.001)
	(162.907,36.291)
	(162.205,32.440)
	(161.486,28.444)
	(160.752,24.298)
	(160.000,20.000)

\path(120,20)	(119.565,23.675)
	(119.419,26.411)
	(120.000,30.000)

\path(120,30)	(121.405,33.049)
	(123.485,36.321)
	(126.322,38.933)
	(130.000,40.000)

\path(130,40)	(133.678,38.933)
	(136.515,36.322)
	(138.595,33.049)
	(140.000,30.000)

\path(140,30)	(140.581,26.411)
	(140.435,23.675)
	(140.000,20.000)

\path(180,20)	(180.579,23.083)
	(181.142,25.733)
	(182.278,29.910)
	(183.525,32.882)
	(185.000,35.000)

\path(185,35)	(189.957,39.476)
	(192.968,41.606)
	(196.265,43.658)
	(199.801,45.624)
	(203.523,47.499)
	(207.382,49.278)
	(211.328,50.955)
	(215.311,52.524)
	(219.280,53.979)
	(223.185,55.315)
	(226.977,56.525)
	(230.604,57.605)
	(234.017,58.547)
	(237.166,59.348)
	(240.000,60.000)

\path(240,60)	(242.546,60.503)
	(245.250,60.939)
	(248.138,61.308)
	(251.239,61.609)
	(254.579,61.844)
	(258.187,62.012)
	(262.089,62.112)
	(266.313,62.146)
	(270.887,62.112)
	(275.838,62.012)
	(278.463,61.936)
	(281.193,61.844)
	(284.031,61.735)
	(286.981,61.609)
	(290.045,61.467)
	(293.228,61.308)
	(296.532,61.132)
	(299.962,60.939)
	(303.520,60.729)
	(307.210,60.503)
	(311.035,60.260)
	(315.000,60.000)

\path(325,60)	(328.320,58.332)
	(331.396,56.755)
	(334.238,55.260)
	(336.859,53.840)
	(341.479,51.196)
	(345.344,48.766)
	(348.540,46.489)
	(351.156,44.309)
	(355.000,40.000)

\path(355,40)	(356.581,37.102)
	(357.863,33.115)
	(358.429,30.567)
	(358.964,27.570)
	(359.483,24.068)
	(360.000,20.000)

\path(400,20)	(400.000,23.034)
	(400.000,25.648)
	(400.000,29.797)
	(400.000,32.797)
	(400.000,35.000)

\path(400,35)	(400.000,37.203)
	(400.000,40.203)
	(400.000,44.352)
	(400.000,46.966)
	(400.000,50.000)

\path(400,65)	(399.296,69.073)
	(398.645,72.579)
	(398.030,75.577)
	(397.438,78.126)
	(396.262,82.110)
	(395.000,85.000)

\path(395,85)	(393.052,88.466)
	(390.767,92.245)
	(388.173,96.280)
	(385.301,100.511)
	(382.181,104.878)
	(378.842,109.322)
	(375.314,113.784)
	(371.628,118.206)
	(367.813,122.527)
	(363.898,126.689)
	(359.914,130.633)
	(355.891,134.298)
	(351.858,137.627)
	(347.845,140.560)
	(343.883,143.037)
	(340.000,145.000)

\path(340,145)	(336.330,146.258)
	(332.049,147.139)
	(327.362,147.652)
	(322.474,147.805)
	(317.589,147.607)
	(312.911,147.068)
	(308.647,146.196)
	(305.000,145.000)

\path(305,145)	(300.014,142.191)
	(297.388,140.110)
	(294.576,137.479)
	(291.504,134.226)
	(288.100,130.277)
	(284.290,125.560)
	(282.209,122.890)
	(280.000,120.000)

\path(260,100)	(256.526,96.022)
	(253.857,93.185)
	(250.000,90.000)

\path(250,90)	(245.180,87.595)
	(242.308,86.433)
	(239.191,85.301)
	(235.878,84.202)
	(232.415,83.138)
	(228.849,82.113)
	(225.229,81.128)
	(221.602,80.185)
	(218.015,79.288)
	(214.516,78.439)
	(211.152,77.640)
	(207.971,76.895)
	(205.021,76.204)
	(202.348,75.572)
	(200.000,75.000)

\path(200,75)	(196.347,74.152)
	(191.353,73.113)
	(188.170,72.487)
	(184.433,71.769)
	(180.067,70.945)
	(177.626,70.488)
	(175.000,70.000)

\path(160,65)	(155.812,64.754)
	(151.918,64.508)
	(148.303,64.260)
	(144.951,64.008)
	(141.849,63.749)
	(138.982,63.484)
	(136.335,63.208)
	(133.894,62.922)
	(129.568,62.307)
	(125.889,61.625)
	(122.739,60.861)
	(120.000,60.000)

\path(120,60)	(115.648,57.990)
	(110.741,55.024)
	(108.137,53.283)
	(105.463,51.422)
	(102.743,49.481)
	(100.000,47.500)
	(97.257,45.519)
	(94.537,43.578)
	(91.863,41.717)
	(89.259,39.976)
	(84.352,37.010)
	(80.000,35.000)

\path(80,35)	(77.350,34.113)
	(74.523,33.279)
	(71.491,32.492)
	(68.223,31.751)
	(64.691,31.050)
	(60.864,30.386)
	(56.715,29.757)
	(52.214,29.157)
	(47.331,28.584)
	(44.737,28.306)
	(42.036,28.033)
	(39.226,27.765)
	(36.302,27.501)
	(33.260,27.241)
	(30.098,26.985)
	(26.811,26.731)
	(23.395,26.480)
	(19.847,26.231)
	(16.164,25.984)
	(12.341,25.737)
	(8.375,25.492)
	(4.263,25.246)
	(0.000,25.000)

\put(210,30){\makebox(0,0)[lb]{
\raisebox{0pt}[0pt][0pt]{\shortstack[l]{{\rm\scriptsize upper half plane}}}}}
\put(210,0){\makebox(0,0)[lb]{
\raisebox{0pt}[0pt][0pt]{\shortstack[l]{{\rm\scriptsize lower half plane}}}}}
\end{picture}
}}
\def\upOV
{{
\begin{picture}(20,35)(0,-10)
\path(15,5)(20,0)
\path(0,0)(20,20)
\path(15.757,12.929)(20.000,20.000)(12.929,15.757)
\path(7.071,15.757)(0.000,20.000)(4.243,12.929)
\path(0,20)(5,15)
\end{picture}
}}
\def\upSI
{{
\begin{picture}(20,35)(0,-10)
\put(10,10){\blacken\ellipse{6}{6}}
\put(10,10){\ellipse{6}{6}}
\path(0,0)(20,20)
\path(15.757,12.929)(20.000,20.000)(12.929,15.757)
\path(20,0)(0,20)
\path(7.071,15.757)(0.000,20.000)(4.243,12.929)
\end{picture}
}}
\def\upUN
{{
\begin{picture}(20,35)(0,-10)
\path(5,5)(0,0)
\path(15.757,12.929)(20.000,20.000)(12.929,15.757)
\path(20,20)(15,15)
\path(7.071,15.757)(0.000,20.000)(4.243,12.929)
\path(0,20)(20,0)
\end{picture}
}}
\begin{document}
\title{Polynomial Invariants are Polynomial}

\author{Dror Bar-Natan}
\address{Department of Mathematics\\
	Harvard University\\
	Cambridge, MA 02138}
\curraddr{Institute of Mathematics\\
        The Hebrew University\\
        Giv'at-Ram, Jerusalem 91904\\
        Israel}
\email{drorbn@math.huji.ac.il}

\if\MRL y
  \thanks{Received January 26, 1995.}
\fi
\thanks{This work was supported by NSF grant DMS-92-03382.}
\if\MRL n
  \thanks{Published in {\em Mathematical Research Letters}, {\bf 2}
    (1995) 239--246, and is available electronically at
    {\tt http://www.ma.huji.ac.il/$\sim$drorbn}, at \newline
    {\tt file://ftp.ma.huji.ac.il/drorbn}, and at {\tt
    http://xxx.lanl.gov/abs/q-alg/9606025}.
}
\fi

\if\MRL n
  \date{This edition: Jun.~30,~1996; \ \ First edition: December 30, 1994.}
\fi

\maketitle

\begin{abstract}
We show that (as conjectured by Lin and Wang) when a Vassiliev invariant
of type $m$ is evaluated on a knot projection having $n$ crossings, the
result is bounded by a constant times $n^m$. Thus the well known
analogy between Vassiliev invariants and polynomials justifies (well,
at least {\em explains}) the odd title of this note.
\end{abstract}

\if\MRL n
  \tableofcontents
\fi
\section{Introduction}

Let $V$ be a fixed Vassiliev knot invariant of type $m$ with values in
some normed vector space (see e.g.~\cite{Bar-Natan:Vassiliev,
Birman:Bulletin, BirmanLin:Vassiliev, Gusarov:New,
Gusarov:nEquivalence, Kontsevich:Vassiliev, Vassiliev:CohKnot,
Vassiliev:Book}). The purpose of this note is to prove the following
three theorems:

\begin{theorem} \lbl{computational}
If $K$ is a knot with a large number $n$ of crossings (in some planar
projection), then $V(K)$ can be computed (in terms of $V$ of finitely
many fixed knots) in $O(n^m)$ computational steps.
\end{theorem}

\begin{theorem} \lbl{special}
If $K$ is a knot with a large number $n$ of crossings (in some planar
projection), then $V(K)$ is bounded by $Cn^m$ for some fixed constant
$C$.
\end{theorem}

\begin{theorem} \lbl{general}
If $K$ is a (singular) knot with $k$ double points and a large number $n$
of crossings (in some planar projection), then $V(K)$ is bounded by
$C_kn^{m-k}$ for some fixed constants $C_k$.
\end{theorem}

We will only prove theorem~\ref{general}. Theorem~\ref{special} follows
from theorem~\ref{general} by setting $k=0$, and
theorem~\ref{computational} can be proven by making all the steps of our
proof effective.

Theorem~\ref{special} was stated as a conjecture
in~\cite{LinWang:IntegralGeometry}, where Lin and Wang commented that it can
be interpreted as saying that polynomial invariants grow polynomially.
Simply recall that in~\cite{Bar-Natan:Vassiliev} an analogy was made
between Vassiliev invariants and polynomials.

\begin{remark} With little additional effort one can generalize the results
of this note to links, tangles, etc.
\end{remark}

\begin{remark} It is rather easy to show that $V(K)$ is bounded by a
polynomial of degree $2m$ or $3m$ in the number of crossings $n$, and
that it is computable in a (high-degree) polynomial time, as stated
in~\cite{Kontsevich:Vassiliev}. For example, one can use the
combinatorial formulas for a universal Vassiliev invariant in terms of
a Drinfel'd associator to find such bounds, or one may argue along the
same lines of this paper but with a little less care about the bounds
in Lemma~\ref{MainLemma}. I found the proof of the much more pleasing
degree $m$ bound for a type $m$ invariant to be somewhat trickier than
expected, as presented in this note.
\end{remark}

\if\MRL n
  \subsection{Acknowledgement}
  \Acknowledgement
\fi
\def\upov{\silentEEPIC{\upOV}{0.35}}     
\def\upun{\silentEEPIC{\upUN}{0.35}}     
\def\upsi{\silentEEPIC{\upSI}{0.35}}     

\section{The method of proof}

For a technical reason (see remark~\ref{WhyLinear}), we prefer to work
with knots parametrized by a parameter $s\in\bold R$ (rather than $s\in
S^1$) and extending from the point $(-\infty,0,-\infty)$ to the point
$(\infty,0,\infty)$ (in some appropriate compactification of $\bold
R^3$). If no double points are allowed, this theory of knots is
equivalent to the usual theory of knots parametrized by a circle.

For any chord diagram $D$ of degree at most $m$ choose a singular knot
$K_D$ representing it (see e.g.~\cite{Bar-Natan:Vassiliev}), and fix
these representatives once and for all. It is well known (see
e.g.~\cite{Birman:Bulletin, BirmanLin:Vassiliev, Vassiliev:CohKnot,
Vassiliev:Book}) that $V$ is determined by its type $m$ and its values
on all the $K_D$'s. This is proven for singular knots with $k$ double
points (``$k$-singular knots'') by downward induction on $k$: If $K$ is
a singular knot with $k>m$ double points then $V(K)=0$ by the defining
property of Vassiliev invariants. And if we know $V(L)$ for every
$(k+1)$-singular knot $L$ (for $k<=m$), we can compute $V(K)$ for a
$k$-singular knot $K$ whose underlying chord diagram is $D$ in terms of
the $V(L)$'s and $V(K_D)$. Simply connect $K$ to $K_D$ via a path
$K(t)$ of singular knots that have exactly $k$ double points at all
times, with the exception of finitely many times where $k+1$ double
points may occur. The usual rule for knots with double points
($V(\upsi)=V(\upov)-V(\upun)$) and telescopic summation now show
that $V(K)-V(K_D)$ is a signed sum of the values of $V$ on the
$(k+1)$-singular knots seen at the exceptional times.

As theorem~\ref{general} is trivial for $k>m$, the above paragraph suggests
that theorem~\ref{general} can be proven for arbitrary $k$ by downward
induction on $k$. Clearly, the induction step (going from $k+1$ to $k$)
follows from the following Lemma:

\begin{Lemma} \lbl{MainLemma}
Let $K_0$ and $K_1$ be two $k$-singular knot projections with
at most $n$ crossings, having the same underlying degree $k$ chord
diagram. Then there exists a path of singular knot projections
connecting $K_0$ and $K_1$, along which there are only finitely many times
in which the number of singular points grows to $k+1$, and so that if
you denote the $(k+1)$-singular knot projections that you see along the
way by $L_1,\ldots,L_p$, then:
\begin{enumerate}
\item $p$ is bounded by a linear function of $n$, whose slope $a_k$ depends
only
  on k.
\item The number of crossings in each of the $L_i$'s is bounded by a linear
   function of $n$, whose slope $b_k$ depends only on $k$.
\end{enumerate}
\end{Lemma}

Indeed, let $K$ be a $k$-singular knot projection with (a large number)
$n$ of crossings and $D$ be its underlying chord diagram, and set
$K_0=K$ and $K_1=K_D$. There is no loss of generality in assuming that
$n$ is larger than the number of crossings in the {\em fixed} knot
$K_D$. Pick a path as in the Lemma, and then by the induction hypothesis
\[ |V(L_i)|<C_{k+1}(b_k n)^{m-k-1},\qquad (1\leq i\leq p). \]
Thus by the discussion in the proceeding paragraph and the bound on $p$,
\[ |V(K)| \leq |V(K_D)|+\sum_{i=1}^p |V(L_i)| \leq
  |V(K_D)|+a_k b_k^{m-k-1} C_{k+1} n^{m-k},
\]
and as there are only finitely many {\em fixed} $K_D$'s to consider, we
find that there is a single constant $C_k$ for which
\[ |V(K)| < C_k n^{m-k} \]
for all $k$-singular knots $K$ with (a large number) $n$ of crossings. \qed
\section{Proof of Lemma 1}

\subsection{A reduction to SubLemmas.}
Let us start with some relevant definitions and SubLemmas.

\begin{definition} We will say that a presentation of a singular knot
$K$ (that is, an appropriate immersion
$K=(K_x,K_y,K_z):\bold R\to\bold R^3$) is {\em
almost monotone} if it satisfies $K_z(s)=s$ for all $s\in\bold R$
except in small neighborhoods of the double points. Notice that
$K$ visits each double point twice, once for a small value of the
parameter $s$ and once for a larger value of $s$. We also require that
$K_z(s)=s$ near those smaller $s$'s, and that near the larger
values of $s$ the knot simply makes a `vertical dive' to meet the
lower strand at the double point, and then climbs vertically up.
Finally, we require that the projection of $K$ to the $xy$-plane
will fall entirely in the upper half plane $\{y>0\}$, except perhaps the
projections of small neighborhoods of some of the double points, which
are allowed to extend just a bit into the lower half plane $\{y\leq 0\}$.
We say that the double points whose projections are in the lower half
plane are {\em exposed}, and if all double points are exposed, we say
that the (almost monotone) presentation $K$ is {\em fully exposed}. See
e.g.\ figures~\ref{VertView} and~\ref{SideView}.
\end{definition}

\begin{figure}[htpb]
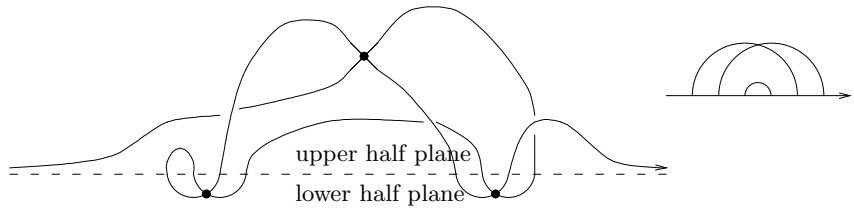

\[ \EEPIC{\VertView}{0.55} \]
\caption{The vertical projection of some almost monotone immersion
  having 3 double points, two of which exposed, and 3 additional crossings.
  On the right is the corresponding chord diagram.
}
\lbl{VertView} \end{figure}

\begin{figure}[htpb]
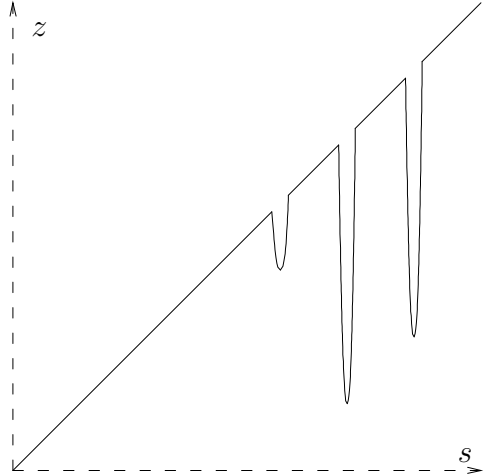

\[ \EEPIC{\SideView}{0.7} \]
\caption{$K_z$ as a function of $s$ for the immersion in
  figure~\ref{VertView}.
}
\lbl{SideView} \end{figure}

\begin{remark} The notion of ``a fully exposed presentation'' is the
key to the proof of Lemma~\ref{MainLemma}. Indeed, within the proof of
SubLemma~\ref{FullyExposed} below, we show that if two fully exposed
presentations have the same underlying chord diagram and their
(exposed) double points (which are in a $1-1$ correspondence) are
embedded in the same way, then the corresponding two singular knots are
the same. In the three SubLemmas below we simply show that any singular
knot presentation can be connected to a fully exposed one by a path
which satisfies the conditions of Lemma~\ref{MainLemma}.
\end{remark}

\begin{SubLemma} \lbl{ToAlmostMonotone}
If a $k$-singular knot $K$ has $n\gg k$ crossings (in some projection), it
can be transformed to an almost monotone knot (having the same chord
diagram, of course) by a path of singular knots satisfying conditions
(1) and (2) of Lemma~\ref{MainLemma}.
\end{SubLemma}

\begin{SubLemma} \lbl{ToFullyExposed}
If a $k$-singular knot presentation $K$ has $n\gg k$ crossings and is
almost monotone, it can be transformed to a fully exposed presentation
(having the same chord diagram, of course) by a path of singular knots
satisfying conditions (1) and (2) of Lemma~\ref{MainLemma}.
\end{SubLemma}

\begin{SubLemma} \lbl{FullyExposed}
If $K_0$ and $K_1$ of Lemma~\ref{MainLemma} have fully exposed
presentations, the conclusion of that Lemma holds.
\end{SubLemma}

Clearly, SubLemmas~\ref{ToAlmostMonotone}--\ref{FullyExposed} imply
Lemma~\ref{MainLemma}. Simply start from $K_0$ and $K_1$, transform them to
be fully exposed using SubLemma~\ref{ToAlmostMonotone} and then
SubLemma~\ref{ToFullyExposed}, and then use SubLemma~\ref{FullyExposed} to
connect the resulting two fully exposed presentations. \qed

\begin{remark} \lbl{WhyLinear} The equality $K_z(s)=s$ in the
definition of almost monotone knots is the reason why it is technically
slightly easier to work with knots parametrized by an infinite line. On
a circle, we'd have to choose some special point where $K_z$ can dive
down so as it can then rise back in a gradual way. Such a point will
have to be given a special treatment, similar to that of the double
points in SubLemma~\ref{ToFullyExposed}, creating some extra mess that
we happily avoid. When dealing with links (as we don't), this extra
mess seems to be unavoidable.
\end{remark}

\subsection{Proofs of the SubLemmas.}

\begin{pf*}{Proof of SubLemma~\ref{ToAlmostMonotone}}
Simply deform $K_z$ to satisfy $K_z(s)=s$ away from the double points
while keeping the projection of $K$ to the $xy$-plane in place. Along
the way you pick some extra double points for crossings that originally
were `the wrong way' (and there are at most $n$ of these), but you
never increase the total number of crossings, so (1) and (2) of
Lemma~\ref{MainLemma} hold. Then (if you're not too tired), do some
cosmetics near the double points to have the strands bounce down and up
as they should.
\end{pf*}

To prove SubLemma~\ref{ToFullyExposed}, we first need

{
  \def\theSubLemma{\ref{ToFullyExposed}}
  \begin{SubSubLemma} \lbl{PathsExist} Let $\pi K$ be the planar
  projection of a $k$-singular knot presentation as in
  SubLemma~\ref{ToFullyExposed} (it is a planar graph with
  $k+n$ vertices and $2(k+n)+1$ edges). There exists disjoint simple
  paths $\gamma_i$ (called `exposing paths')  connecting the projections
  of the un-exposed double points of $K$ to points in the lower half
  plane, so that:
  \begin{itemize}
  \item The $\gamma_i$'s miss all the vertices of $\pi K$.
  \item  The total number of intersection points between the $\gamma_i$'s
    and the edges of $\pi K$ is at most $k(2(k+n)+1)$.
  \end{itemize}
  Check figure~\ref{OddChoice} for an example.
  \end{SubSubLemma}
}

\begin{figure}[htpb]
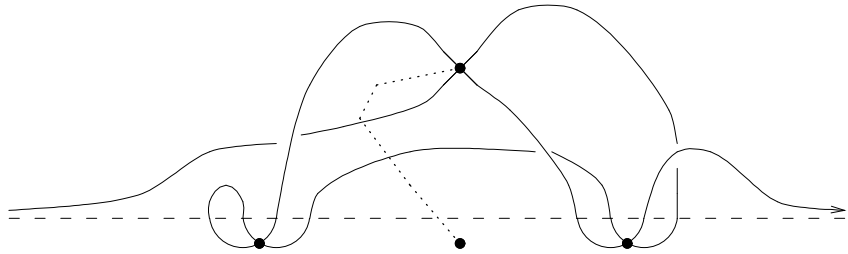

\[ \EEPIC{\OddChoice}{0.7} \]
\caption{A (somewhat odd) choice for an exposing path for the only
  un-exposed double point in the knot projection of
  figure~\ref{VertView}.
}
\lbl{OddChoice} \end{figure}

\begin{pf} Start with arbitrary paths that miss the vertices of $\pi K$
and connect  the projections of the un-exposed double points of $K$ to
points in the lower half plane. If any of these paths intersects any of
the edges of $\pi K$ more than once, at least one of these intersection
points can be eliminated by traveling from one to the other along a
`shorter' path that follows the relevant edge:
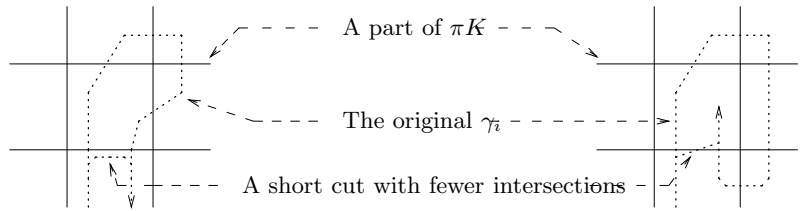
\[ \EEPIC{\Cutting}{0.6} \]
Doing as much as we can of that, we get a collection of paths, each of
which intersecting each edge of $\pi K$ at most once, to a total of at most
$k(2(k+n)+1)$ intersections. But we may have created lots of intersections
between the different paths and lots of self intersections. Eliminate these
by moves like
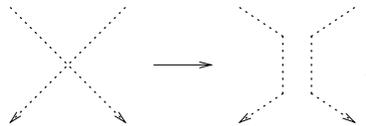
\[ \EEPIC{\Moves}{0.6}, \]
and by throwing out closed loops when these are created as a move like
above is applied to a self-intersection.
\end{pf}

\begin{pf*}{Proof of SubLemma~\ref{ToFullyExposed}}
Choose exposing paths $\gamma_i$ as in SubSubLemma \ref{PathsExist},
and pull the un-exposed double points (and small neighborhoods thereof)
along them:
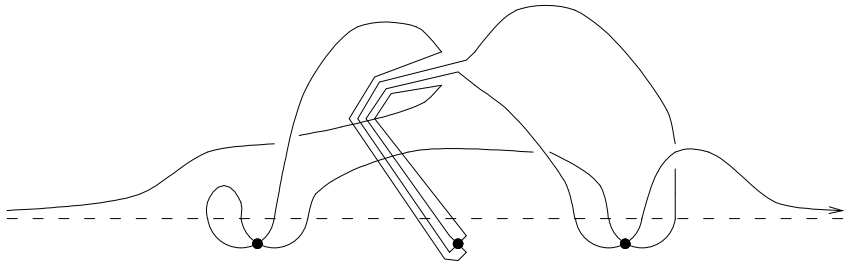
\[ \EEPIC{\Pulling}{0.7} \]
The bound supplied by SubSubLemma~\ref{PathsExist} on the number of
intersections between the $\gamma_i$'s and the projection of $K$ shows
that conditions (1) and (2) of Lemma~\ref{MainLemma} hold.
\end{pf*}

\begin{pf*}{Proof of SubLemma~\ref{FullyExposed}} The fact that $K_0$
and $K_1$ have the same underlying chord diagram implies that there is
a natural correspondence between their $k$ double points, and between
the $4k$ strands emanating from these $k$ double points in each of
them.  Ensure that these $4k$ strands on $K_1$ enter the upper half
plane in the same places as the corresponding ones for $K_0$. This can
be done by permuting and rotating the $k$ double points of $K_1$, at a
cost (in the sense of conditions (1) and (2) of Lemma~\ref{MainLemma})
proportional to $k^2$, not even linear in $n$. The new $K_1$ is now the
{\em same} knot as $K_0$. Indeed, we've just arranged things so that
the restrictions $\pi K_0'$ and $\pi K_1'$ of their projections to the
upper half plane are homotopically equivalent modulo the boundary.  But
both knots are almost monotone, and thus we can lift any homotopy that
takes $\pi K_0'$ to $\pi K_1'$ to an isotopy taking $K_0$ to $K_1$.
\end{pf*}
\ifx\undefined\bysame
        \newcommand{\bysame}{\leavevmode\hbox to3em{\hrulefill}\,}
\fi

\end{document}